\documentclass[a4paper,10pt]{amsart}
\usepackage{graphicx} 
\usepackage{amssymb}

\usepackage{tikz, tikz-3dplot} 

\usepackage[backref=page]{hyperref}

\numberwithin{equation}{section}

\usepackage{color}

\newtheorem{thm}{Theorem}[section]
\newtheorem{prop}[thm]{Proposition}
\newtheorem{lem}[thm]{Lemma}
\newtheorem{cor}[thm]{Corollary}
\newtheorem{rem}[thm]{Remark}

\newtheorem{dfn}[thm]{Definition}

\newcommand{\bC}{{\mathbb C}}

\newcommand{\bQ}{{\mathbb Q}}
\newcommand{\bZ}{{\mathbb Z}}
\newcommand{\bN}{{\mathbb N}}

\newcommand{\mb}{\mathbf}
\newcommand{\mr}{\mathrm}

\title[GKM Structures on Quiver Grassmannians]{Skeletal Torus Actions and GKM Structures on Quiver Grassmannians of String Representations}
\author{Alexander Pütz}
\address{A. P\"utz:\newline
Faculty of Mathematics\\
Bielefeld University\\
PO Box 100 131\\\newline
 D-33501 Bielefeld\\
Germany}
\email{alexander.puetz@math.uni-bielefeld.de}

\begin{document}
\begin{abstract}
Quiver Grassmannians of equioriented type $\texttt{A}$ and nilpotent equioriented type $\tilde{\texttt{A}}$ quiver representations are GKM-varieties. In particular, they have a cellular decomposition and admit a torus action with finitely many fixed points and one-dimensional orbits (i.e. skeletal action). We examine the case of string representations and provide a classification of all corresponding quiver Grassmannians with a GKM-variety structure.
\end{abstract}
\maketitle
\section{Introduction}
Quiver Grassmannians arise naturally as generalization of classical Grassmannians and flag varieties. Given a quiver representation, they parameterize all subrepresentations of a fixed dimension. For example, (degenerate) flag varieties of type $\texttt{A}$ are quiver Grassmannians for equioriented quivers of type $\texttt{A}$ \cite{CFR12}. A projective variety is a GKM-variety if an algebraic torus acts on it with finitely many fixed points and finitely many one-dimensional orbits and the rational cohomology vanishes in odd degrees.  In \cite{CFR13a} the quiver Grassmannian describing the Feigin degeneration of the type $\texttt{A}$ flag variety is equipped with a GKM-variety structure. It was discovered in \cite{LP23a} that this generalizes to quiver Grassmannians for nilpotent representations of the equioriented cycle. This setting includes all representations of equioriented quivers of type $\texttt{A}$. Hence it is natural to ask: Do quiver Grassmannians for representations of other quivers admit a GKM-variety structure? 

In this paper, we provide a full answer for quiver Grassmannians of string representations (Theorem~\ref{thm:GKM-structure} and Lemma~\ref{lem:no-GKM-structure}). Moreover, we include an example suggesting that there exists a class of tree representations with GKM-variety structures on their quiver Grassmannians (Appendix~~\ref{app:example}). However, this case requires new methods and should be investigated independently.

One way to prove vanishing of odd-degree cohomology is to establish a cellular decomposition. For quiver Grassmannians this goal is tackled in \cite{CEFR2018}. They prove that for representations of Dynkin quivers, every quiver Grassmannian admits a cellular decomposition and conjecture the same for quivers which are orientations of the simply laced Dynkin diagrams of affine type. In this paper we do not restrict to the (affine) Dynkin type. Note that cellular decompositions can not exist for arbitrary representations of wild quivers since every projective variety has a realization as quiver Grassmannian for a suitable representation of any wild quiver \cite{Ringel2018}. 

Observe that for quiver Grassmannians it is sufficient to show that the fixed point set of $\mathbb{C}^*$ (acting as subgroup of the higher rank torus) is finite. Then both fixed point sets coincide automatically, since their cardinality equals the Euler characteristic of the variety. 
$\mathbb{C}^*$-actions on quiver representations and their Grassmannians were studied in \cite{Cerulli2011,Haupt2012}. They prove that string, tree and band representations of arbitrary quivers admit a $\mathbb{C}^*$-action such that the number of fixed points in any quiver Grassmannian associated to these representations is finite. Additionally, the fixed points are described by so called coordinate subrepresentations, i.e. their vector spaces are spanned by subsets of the appropriate cardinality in the bases of the representation defining the quiver Grassmannian. However, the actions in \cite{Haupt2012} do not extend to the quiver Grassmannian in general (see Lemma~\ref{lem:all-nice-gradings-are-constructible}).

In order to examine the structure of the one-dimensional torus orbits we want a cellular decomposition into the attracting sets of the $\mathbb{C}^*$-fixed points in the sense of \cite{Birula1973}. Then it is possible to compute the torus orbits in the attracting sets explicitly. 
Furthermore, the constructions (for the $\mathbb{C}^*$-action and cellular decomposition) presented in this paper and in \cite[Section~5.1]{LP23a} should also apply to the setting of certain tree representations. For more details on this,, see \cite[Section~8.2]{LP23a}.

Quiver Grassmannians of type $\texttt{A}$ are closely related to the type $\texttt{A}$ flag variety (see \cite{CFR12,CFFFR17}). Hence, it is natural to ask if the description of full exceptional collections in the bounded derived category of coherent sheaves on generalized flag varieties $G/P$ from \cite{SvdK24} extends to quiver Grassmannians. A first step in this direction is a good understanding of the GKM-cohomology resp. $T$-equivariant $K$-theory because the elements of the full exceptional collection are parameterized by the $T$-fixed points (see \cite[Section~1.1]{SvdK24}). For equioriented type $\texttt{A}$ and $\tilde{\texttt{A}}$ quiver Grassmannians, the GKM-cohomology was examined in \cite{LP23a,LP23b}. 
This is a first step towards a description of their full exceptional collections. A natural case to start the study of full exceptional collections on quiver Grassmannians are linear degenerations of flag varieties of type $\texttt{A}$ (see \cite{CFFFR17} for their definition).


This paper is structured as follows. In Section~\ref{sec:Background} we recall some definitions and results concerning quiver Grassmannians and GKM theory. Section~\ref{sec:torus-actions} is devoted to the study of torus actions on quiver Grassmannians for representations that are the direct sum of strings and trees. We prove that all gradings whose induced actions extend to the quiver Grassmannian are constructible from an initial weight of every indecomposable direct summand and weights of the edges in the quiver (Lemma~\ref{lem:all-nice-gradings-are-constructible} and Remark~\ref{rem:constructing-constructible-gradings}). In Section~\ref{sec:GKM-Quiver-Grass} we equip quiver Grassmannians for equioriented string representations of arbitrary quivers with a GKM-variety structure (Theorem~\ref{thm:GKM-structure}). This is the main result of the present paper. Its proof is divided into the construction of a cellular decomposition (Theorem~\ref{thm:Affine-Paving}), the description of torus fixed points (Lemma~\ref{lem:finitely-many-fixed-points}) 
and one-dimensional torus orbits (Theorem~\ref{thm:one-dim-T-orbits}). Moreover, we prove that in general a quiver Grassmannian for an arbitrary string representation can not have a GKM-variety structure if it is not equioriented (Lemma~\ref{lem:no-GKM-structure}).  In Section~\ref{sec:momnent-graph} we apply the GKM-variety structure described in Section~\ref{sec:GKM-Quiver-Grass}.  The moment graph for the torus action from Theorem~\ref{thm:GKM-structure} has a combinatorial description via certain cut and paste moves on the coefficient quivers of the fixed points, called fundamental mutations (Theorem~\ref{thm:moment-graph}). 
We compute a combinatorial basis for the torus equivariant cohomology of the quiver Grassmannians from Theorem~\ref{thm:GKM-structure} in terms of Knutson-Tao classes (Theorem~\ref{lem:KT-basis}).
The tangent space of the quiver Grassmannian at a torus fixed point has a basis given by the mutations of the coefficient quiver corresponding to the fixed point (Theorem~\ref{thm:tangent-space-from-mutations}).
Finally, in Appendix~\ref{app:example} we give an example of a quiver Grassmannian with GKM-variety structure for a representation which is not a string.
\section{Background}\label{sec:Background}
\subsection{Quiver representations and Grassmannians}
For more detail beyond this brief summary, we refer to \cite{Cerulli2011,Kirillov2016,Schiffler2014}. Let $Q$ be a quiver, consisting of a set of vertices $Q_0$ and a set of arrows $Q_1$ between the vertices. The source and target of an arrow $a \in Q_1$ are denoted by $s_a,t_a \in Q_0$. A $Q$-representation $M$ consists of a tuple of $\bC$-vector spaces $M^{(i)}$ for $i \in Q_0$ and tuple of linear maps $M_a : M^{(i)} \to M^{(j)}$ for $a \in Q_1$ with $i =s_a$ and $j=t_a$. By $\mr{rep}_\bC(Q)$ we denote the category of finite dimensional $Q$-representations. Morphisms between two objects $N$ and $M$ are tuples of linear maps $\varphi_i : N^{(i)} \to M^{(i)}$ for $i \in Q_0$ such that 
\( \varphi_j \circ N_a = M_a \circ \varphi_i\) holds for all $a \in Q_1$ with $i =s_a$ and $j=t_a$. The set of all such morphisms is $\mr{Hom}_Q(N,M)$. By abuse of notation we write $M \in \mr{rep}_\bC(Q)$ for the objects of the category. The path algebra of the quiver $Q$ is denoted by $\bC Q$. Its underlying vector space is spanned by all paths in the quiver and multiplication is given by concatenation of paths. 
\begin{dfn}
For $M \in \mr{rep}_\bC(Q)$ and $\mb{e} \in \bN^{Q_0}$, the \textbf{quiver Grassmannian} $\mr{Gr}_\mb{e}(M)$ is the closed subvariety in 
$\prod_{i\in Q_0}{\rm Gr}_{e_i}(M^{(i)})$ parameterizing all subrepresentations $U \subseteq M$ such that $\dim_\bC U^{(i)}=e_i$ for $i \in Q_0$.
\end{dfn}
\subsection{Tree representations}

A {quiver morphism} $F: S \to Q$ between two finite quivers $S$ and $Q$ is a pair of maps $F_0: S_0 \to Q_0$ and $F_1: S_1 \to Q_1$ such that for each arrow $(a : i \to j) \in S_1$ its image $F_1(a)$ is the arrow $(b: F_0(i) \to F_0(j)) \in Q_1$.
We call a quiver morphism $F: S \to Q$ \textbf{winding} if in $S$ there are no subquivers of the form
\begin{center}
\begin{tikzpicture}[scale=.8]
   
	\draw[fill=black] (-1,0) circle(.08);
	\draw[arrows={-angle 90}, shorten >=2.5]  (-1,0) -- (0,0); \node at (-0.5,0.2) {$a$};
	\draw[fill=black] (0,0) circle(.08);
	\draw[arrows={-angle 90}, shorten >=2.5]  (1,0) -- (0,0); \node at (0.5,0.26) {$b$};
	\draw[fill=black] (1,0) circle(.08);
	
\end{tikzpicture} $\quad \quad \mr{or}\quad \quad$
\begin{tikzpicture}[scale=.8]
   
	\draw[fill=black] (-1,0) circle(.08);
	\draw[arrows={-angle 90}, shorten >=2.5]  (0,0) -- (-1,0); \node at (-0.5,0.2) {$a$};
	\draw[fill=black] (0,0) circle(.08);
	\draw[arrows={-angle 90}, shorten >=2.5]  (0,0) -- (1,0); \node at (0.5,0.26) {$b$};
	\draw[fill=black] (1,0) circle(.08);
	
\end{tikzpicture} 
\end{center} 
with $F(a) = F(b)$. A winding $F$ induces a push-down functor $F_* : \mr{rep}_\bC (S) \to \mr{rep}_\bC (Q)$ defined by
\[ (F_*V)^{(i)} := \bigoplus_{j \in F^{-1}(i)} V^{(j)} \ \ \mr{for} \ i \in Q_0 \quad \quad (F_*V)_a := \bigoplus_{b \in F^{-1}(a)} V_b \quad \ \ \mr{for} \ a \in Q_1. \]
\begin{dfn}\label{def:tree-rep}
A $Q$-representation of the form $F_*V_S$, where $S$ is a tree and $F$ is a winding, is called \textbf{tree representation}. 
Here $V_S$ we denotes the $S$-representation, with $V_S^{(i)}=\bC$ for all $i \in S_0$ and  $V_{S,a}= \mr{id}_\bC$ for all $a \in S_1$. 
If $S$ is a type $A$ quiver, $F_*V_S$ is called \textbf{string representation}. A direct sum of trees is called \textbf{forest}.
\end{dfn}
The following combinatorial object takes a cental role in this paper.
\begin{dfn}\label{dfn:coeff-quiver}
Let $M \in \mathrm{rep}_\mathbb{C}(Q)$ and let $B=\{b^{(i)}_k\}$ denote a basis of the vector space $V = \oplus_{i \in Q_0} M^{(i)}$. The \textbf{coefficient quiver} $Q(M,B)$ consists of: 
\begin{itemize}
\item[(QM0)] the vertex set  $Q(M,B)_0=B$,
\item[(QM1)] the set of arrows $Q(M,B)_1$, containing $(\tilde{a}: b^{(i)}_k \to b^{(j)}_\ell)$ if and only if $(a:i\to j) \in Q_1$ and the coefficient of $b^{(j)}_\ell$ in $M_a (b^{(i)}_k)$ is non-zero.
\end{itemize}    
\end{dfn}
\begin{rem}
The coefficient quiver of $F_*V_S$ is $S$.
\end{rem}
\begin{rem}\label{rem:different-def-of-trees}
Definition~\ref{def:tree-rep} follows \cite[p.~260]{C-B1989} and differs from other literature, where it is only required that there exists a basis $B$ of $M \in \mr{rep}_\bC(Q)$ such that $Q(M,B)$ is a tree (see \cite[Remark~5]{Ringel1998}). 
\end{rem}

\begin{rem}\label{rem:trees-are-indec}
By \cite[Section~3.5 and Section~4.1]{Gabriel1981}, tree representations as in Definition~\ref{def:tree-rep} are indecomposable. This does not hold for the weaker assumption on the coefficient quiver. In \cite{Cerulli2011} this is overcome by the notion of orientable strings. Their definition is equivalent to Definition~\ref{def:tree-rep}.
\end{rem}
\subsection{Basics of GKM-Theory}\label{sec:Basics-GKM-Theory} 
The notation of this section is based on \cite{LP23b}. The original theory is found in \cite{GKM98}.
Let $X$ be a projective algebraic variety over $\bC$. The action of an algebraic torus $T \cong (\bC^*)^r$ on $X$ is \textbf{skeletal} if the number of $T$-fixed points and the number of one-dimensional $T$-orbits in $X$ is finite. A cocharacter $\chi \in \mathfrak{X}_*(T)$ is called \textbf{generic} for the $T$-action on $X$ if $X^T = X^{\chi(\bC^*)}$. By $H_T^\bullet(X)$ we denote the $T$-equivariant cohomology of $X$ with rational coefficients. 
\begin{dfn}\label{dfn:GKM-variety} The pair $(X,T)$ is a \textbf{GKM-variety} if the $T$-action on $X$ is skeletal and the rational cohomology of $X$ vanishes in odd degrees.
\end{dfn}
For every one-dimensional $T$-orbit $E$ in a projective GKM-variety there exists a $T$-equivariant isomorphism between its closure $\overline{E}$ and $\bC\mathbb{P}^1$. Thus each one-dimensional $T$-orbit connects two distinct $T$-fixed points of $X$.
\begin{dfn}\label{def:moment-graph}Let $(X,T)$ be a GKM-variety, and let $\chi\in \mathfrak{X}_*(T)$ be a generic cocharacter. The corresponding \textbf{moment graph} $\mathcal{G}=\mathcal{G}(X,T, \chi)$ of a GKM-variety is given by the following data:
\begin{itemize}
\item[(MG0)] the $T$-fixed points as vertices, i.e.: $\mathcal{G}_0 = X^T$,
\item[(MG1)] the closures of one-dimensional $T$-orbits $\overline{E} = E \cup \{x,y\}$ as edges in $\mathcal{G}_1$, oriented from $x$ to $y$ if $\lim_{\lambda \to 0}\chi(\lambda).p=x$ for $p \in E$,
\item[(MG2)] every $\overline{E}$ is labeled by a character $\alpha_E \in \mathfrak{X}^*(T)$ describing the $T$-action on $E$. 
\end{itemize}
\end{dfn}
\begin{rem}\label{rem:torus-characters-part-i}
Fixing a generic cocharacter $\chi\in \mathfrak{X}_*(T)$, determines the orientation of the edges in (MG1). The characters in (MG2) are unique up to a sign but this sign plays no role in Theorem~\ref{thm:GKM}. Hence we fix a choice from now on. 
\end{rem}
\begin{rem}\label{rem:torus-characters-part-ii} Let $T$ be a torus of rank $q$ and let $\{\tau_1, \ldots, \tau_q\}$ be a $\bZ$-basis  of its character lattice $\mathfrak{X}^*(T)$. For any torus character $\alpha$ we also write $\alpha$ for its image $\alpha\otimes 1$ in the $\bQ$-vector space $\mathfrak{X}^*(T)\otimes_\bZ \bQ$.
Identify the symmetric algebra of this vector space with $\bQ[T]=\bQ[\tau_1,\dots,\tau_q]$, and the latter with $R:=H_T^\bullet(\textrm{pt})$. 
\end{rem}
One main motivation to study moment graphs is the following result: 
\begin{thm}\label{thm:GKM}(\cite{GKM98}) Let $(X, T)$ be a GKM-variety with moment graph $\mathcal{G}=\mathcal{G}(X,T,\chi)$. Then 
\[
H_T^\bullet(X) \cong \left\{(f_x)\in\bigoplus_{x\in \mathcal{G}_0}R \ \Big| \
\begin{array}{c}
f_{x_E}-f_{y_E}\in \alpha_{E} R\\
 \hbox{ for any }\overline{E}=E\cup\{x_E, y_E\}\in\mathcal{G}_1\end{array}
\right\}.\]  
\end{thm}
\section{Torus actions}\label{sec:torus-actions}
In this section we introduce torus actions on the underlying vector space of a quiver representation and examine under which assumption it extends to the corresponding quiver Grassmannians. So far, all quiver Grassmannians with GKM-variety structure belong to representations of quivers of type \texttt{A} and type $\tilde{\texttt{A}}$ \cite[Theorem~6.6]{LP23a}. In particular, these representations are strings. Hence, it is natural to start with the examination of torus actions on string and tree representations of arbitrary quivers. From now on, let $M \in \mathrm{rep}_\bC (Q)$ be a forest.
\subsection{Constructible gradings}
For $M \in \mathrm{rep}_\bC(Q)$, let $B$ denote a basis of the vector space $V = \oplus_{i \in Q_0} M^{(i)}$. A \textbf{grading} is a map $\mr{wt}: B \to \bZ$. It induces a $\bC^*$-action on $V$ by 
\[ \lambda.b := \lambda^{\mr{wt}(b)}\cdot b \quad \text{for all } b \in B \text{ and all }\lambda \in \bC^*.  \]  
A grading is called \textbf{e-admissible} if its induced $\bC^*$-action extends to $\mr{Gr}_\mb{e}(M)$, for a fixed $\mb{e} \in \bZ^{Q_0}$.  
\begin{rem}
 The action of higher rank tori is obtained via the combination of multiple gradings.
\end{rem}
A grading is called \textbf{constructible} if every arrow $a \in Q_1$ has a fixed weight $\mr{wt}(a) \in \bZ$, i.e. $\mr{wt}(b')=\mr{wt}(b)+\mr{wt}(a)$ for all $b,b' \in B$ such that $M_a b=\mu b'$ for some $\mu \in \bC^*$. 
\begin{prop}\label{prop:strong-admissible-gradings}
    Constructible gradings are admissible for all $\mb{e} \in \bZ^{Q_0}$.
\end{prop}
\begin{proof}
We have to check that $\lambda.U \in \mr{Gr}_\mb{e}(M)$ for all $\lambda \in \bC^*$ and all $U \in \mr{Gr}_\mb{e}(M)$. This holds if and only if $M_a (\lambda.U^{(i)} )\subseteq \lambda.U^{(j)}$ holds for all $a \in Q_1$ with $i = s_a$ and $j=t_a$. Let $m_i := \dim M^{(i)}$ and let $\{b^{(i)}_t : \ell \in[m_i]\}$ denote the basis of $M^{(i)}$, then the basis vectors of $U^{(i)}$ are 
\[ u^{(i)}_k=\sum_{\ell \in [m_i]}\mu^{(i)}_{k,\ell}b^{(i)}_\ell \]
for suitable $\mu^{(i)}_{k,\ell} \in \bC$ and $k \in [e_i]$. Since $M$ consists of windings of trees, for all $a \in Q_1$ with $i = s_a$ and $k \in [m_i]$ there exists $k' \in [m_j]$ for $j=t_a$ and $\nu_{a,k,k'} \in \bC^*$ such that $M_a (  b^{(i)}_k ) = \nu_{a,k,k'} b^{(j)}_{k'}$. Now, the constructibility of the grading implies
\begin{align*} M_a \big(  \lambda.u^{(i)}_k \big) &=  M_a \Big(  \lambda. \sum_{\ell \in [m_i]}\mu^{(i)}_{k,\ell}b^{(i)}_\ell \Big) = M_a \Big( \sum_{\ell \in [m_i]} \mu^{(i)}_{k,\ell}\lambda^{\mr{wt}\big(b^{(i)}_\ell\big)}   b^{(i)}_\ell \Big)\\
&=  \sum_{\ell \in [m_i]} \mu^{(i)}_{k,\ell}\lambda^{\mr{wt}\big(b^{(i)}_\ell\big)} 
 M_a \big(  b^{(i)}_\ell \big) = \sum_{\ell \in [m_i]} \mu^{(i)}_{k,\ell}\lambda^{\mr{wt}\big(b^{(i)}_\ell\big)} 
 \nu_{a,\ell,\ell'}  b^{(j)}_{\ell'} \\
 &= \sum_{\ell \in [m_i]} \mu^{(i)}_{k,\ell} \nu_{a,\ell,\ell'} \lambda^{-\mr{wt}\big(a\big)} \lambda^{\mr{wt}\big(b^{(i)}_{\ell'}\big)} 
   b^{(j)}_{\ell'} = \lambda^{-\mr{wt}\big(a\big)} \sum_{\ell \in [m_i]} \mu^{(i)}_{k,\ell} \lambda.  M_a \big(  b^{(i)}_\ell \big)\\
   &=\lambda^{-\mr{wt}\big(a\big)} \cdot \lambda.M_a \big(  u^{(i)}_k \big).
\end{align*} 
Hence $\lambda.U \in \mr{Gr}_\mb{e}(M)$ is equivalent to $U \in \mr{Gr}_\mb{e}(M)$ for all $\lambda \in \bC^*$.
\end{proof}
\begin{rem}
    We refer to gradings which are admissible for all $\mb{e} \in \bZ^{Q_0}$ as \textbf{strong admissible}. It is natural to ask for other classes of (strong) admissible gradings. In \cite{Haupt2012}, the class of so called nice gradings is introduced. Unfortunately, they are only strong admissible if they are already constructible (see Lemma~\ref{lem:all-nice-gradings-are-constructible}).
\end{rem}
\subsection{Flexible subrepresentations}
\begin{dfn}
 A dimension vector $\mb{e} \in \bZ^{Q_0}$ is called \textbf{M-flexible} if for all $i \in Q_0$, $\dim_\bC  \mr{pr}_i \mr{Gr}_\mb{e}(M)>0$, where $\mr{pr}_i: \mr{Gr}_\mb{e}(M)  \to \mr{Gr}_{e_i}(M^{(i)})$, $U \mapsto U^{(i)}$.   
\end{dfn}
\begin{prop}\label{prop:flexible-reduction}
    If $i \in Q_0$ is not $M$-flexible, 
    then there exist $M' \in \mr{rep}_\bC (Q')$ and $\mb{e}' \in \bZ^{Q_0'}$ for $Q' := Q \setminus \{i\}$ such that $\mr{Gr}_\mb{e}(M) \cong \mr{Gr}_\mb{e'}(M')$.
\end{prop}
\begin{proof}
    Since $U^{(i)}=M^{(i)}$ for all $U \in \mr{Gr}_\mb{e}(M)$, the maps $M_a$ with $t_a = i$ do not influence the variety and can be removed. Let $B$ denote the basis for the vector spaces of $M$. For all paths $(p:i \to j) \in \bC (Q)$, remove the basis vectors in the image of $M_p(M^{(i)})$ from $B$ and decrease $e_j$ by $\dim M_p(M^{(i)})$. Here, $M_p$ denotes the concatenation of the maps $M_a$ along the path $p$.
    
    Let $M'$ have vector spaces spanned by the remaining basis vectors and take the maps of $M$ restricted to them. The maps $M_a$ with $s_a =i$ are removed. It is straight forward to check that there is a pointwise isomorphism between $\mr{Gr}_\mb{e}(M)$ and $\mr{Gr}_\mb{e'}(M')$. 
\end{proof}
\begin{rem}
    It is always possible to repeat this reduction until either $\mb{e'}$ is $M'$-flexible or both are trivial. The latter occurs when the starting quiver Grassmannian was a point or empty.  
\end{rem}
\begin{lem}\label{lem:all-nice-gradings-are-constructible}
    If $\mb{e} \in \bZ^{Q_0}$ is $M$-flexible, every $\mb{e}$-admissible grading is constructible.
\end{lem}
\begin{proof}
Assume that $\mr{wt} : B \to \bZ$ is not constructible, i.e. there exists an arrow $a \in Q_1$ and $b_1,b_2 \in B$ with $M_a b_1 \neq 0$ and $M_a b_2 \neq 0$ such that $\mr{wt}( M_a b_1) - \mr{wt}(b_1) \neq \mr{wt}( M_a b_2) - \mr{wt}(b_2)$. 
At vertex $i =s_a$ consider the subspace spanned by 
\[ u^{(i)}_k=\sum_{\ell \in [m_i]}\mu^{(i)}_{k,\ell}b^{(i)}_\ell \]
such that for exactly one $k_0 \in [e_i]$ we have $\mu^{(i)}_{k_0,\ell} \neq 0$ for $b^{(i)}_\ell=b_1$ and $b^{(i)}_\ell=b_2$. All other $u^{(i)}_k$ satisfy $\mu^{(i)}_{k_0,\ell} = 0$ for $b^{(i)}_\ell=b_2$. Since $\mb{e} \in \bZ^{Q_0}$ is $M$-flexible, this choice exists. It extends to a representation $U \in \mr{Gr}_\mb{e}(M)$ with analogous condition at vertex $j = t_a$, because $M$ consists of windings of trees and $M_a b_1 \neq 0$, $M_a b_2 \neq 0$. By assumption $u^{(i)}_{k_0}$ and $u^{(j)}_{k_0'}$ are not contained in the span of the remaining basis vectors for the vector spaces of $U$. It is computed analogous to the proof of Proposition~\ref{prop:strong-admissible-gradings}, that there exists $\lambda \in \bC^*$ such that $M_a (\lambda.U^{(i)}) \nsubseteq \lambda.U^{(j)}$ and hence $\lambda.U \notin \mr{Gr}_\mb{e}(M)$. Consequently, $\mr{wt} : B \to \bZ$ is not $\mb{e}$-admissible.
\end{proof}
\begin{rem}
In particular, the action of arbitrary rank tori on quiver Grassmannians for forests is obtained from constructible gradings. 
\end{rem}
\begin{rem}\label{rem:constructing-constructible-gradings}
 Constructible gradings are uniquely determined by the edge weights $\mr{wt}(a) \in \bZ$ for all $a \in Q_1$ and initial weights $k_j \in \bZ$ for $j \in [d]$ assigned to a distinguished basis vector $b_j \in B$ for each indecomposable summand of $M$.   
\end{rem}
\begin{cor}\label{cor:faithful-action}
Let $M \in \mathrm{rep}_\bC (Q)$ be a forest and let $\mb{e} \in \bZ^{Q_0}$ be $M$-flexible. The rank of the maximal torus acting faithfully on $\mr{Gr}_\mb{e}(M)$ equals $d+c$, where $d$ is the number of indecomposable summands of $M$ and $c = \# (Q_1\cap\mr{supp}M)$.
\end{cor} 
\section{Quiver Grassmannians with GKM Variety Structure}\label{sec:GKM-Quiver-Grass}
\begin{dfn}
    A representation $M \in \mathrm{rep}_\bC (Q)$ is called \textbf{straight} if there exists a basis $B$ for the vector spaces of $M$ such that all connected components of the coefficient quiver $Q(M,B)$ are equioriented strings. A compatible basis $B$ is called \textbf{straightening basis}. 
\end{dfn}
\begin{rem}
This does not require that $Q$ is equioriented. In particular, loops and cycles are allowed. 
\end{rem}
\begin{dfn}\label{dfn:torus-action}
Let $M \in \mathrm{rep}_\bC (Q)$ be straight with straightening basis $B$. 
By $d$ we denote the number of components in $Q(M,B)$ and $c$ is the number of edges in the $M$-supporting subquiver $\mr{supp}M \subseteq Q$. The action of the torus $T_M := (\bC^{*})^{d+c}$ on the vector space basis is defined inductively: On the start of the $i$-th equioriented strings it acts with $\gamma_i$. Along every edge in $Q(M,B)$ we multiply by $\nu_j$ corresponding to the underlying edge of $Q$ labeled by $j \in [c]$.
\end{dfn}
\begin{rem}\label{rem:T-action-in-type-A} 
By Corollary~\ref{cor:faithful-action}, this is the largest torus acting faithfully. For equioriented quivers of type $A$, we recover the torus from \cite{CFR13a} by setting $\nu_j=1$ for all $j \in [c]$. The action on nilpotent representations of the equioriented cycle $\Delta_n$ as constructed in \cite{LP23a} is recovered with $\nu_j=\gamma_0$ for all $j \in [c]$. \end{rem}
\begin{thm}\label{thm:GKM-structure}
    Let  $M \in \mathrm{rep}_\bC (Q)$ be straight and $\mb{e} \in \bZ^{Q_0}$. The action of $T_M$ as in Definition~\ref{dfn:torus-action} equips $\mr{Gr}_\mb{e}(M)$ with a GKM variety structure.
\end{thm}
The rest of the section is devoted to the proof of this theorem.
\subsection{Torus fixed points}
\begin{lem}\label{lem:finitely-many-fixed-points}
    Let  $M \in \mathrm{rep}_\bC (Q)$ be straight and $\mb{e} \in \bZ^{Q_0}$. The torus $T_M$ as in Definition~\ref{dfn:torus-action} acts on $\mr{Gr}_\mb{e}(M)$ with finitely many fixed points.
\end{lem}
\begin{proof}
By \cite[Theorem~1]{Cerulli2011}, a constructible grading induces a $\bC^*$-action with finitely many fixed points if all elements of $B$ 
have distinct weights. Hence, we have to find a constructible cocharacter of $T_M$ with this property.

Since $M$ is straight, it is the direct sum of equioriented strings $s_1,\dots,s_d$. Let $b_{j,1}$ for $j \in [d]$ denote the starting vertex of the $j$-th string and let $\ell$ be the number of vertices on the longest string. We set $\mr{wt}(a)=1$ for all $a \in Q_1$, and $\mr{wt}(b_{j,1})=(j-1)\ell$ for all $a \in Q_1$. This is sufficient to see that all weights are distinct. Clearly, this describes a cocharacter of $T_M$. 
\end{proof}
\subsection{Affine paving}
Let $X$ be a projective variety and let $\bC^*$ act on $X$ with finitely many fixed points. Let $\{x_1, \ldots, x_m\}$ be the fixed point set, which we denote by $X^{\bC^*}$. This induces a decomposition
\begin{equation}\label{eqn:BBdecomposition}
X=\bigcup_{i\in [m]} W_i, \quad\hbox{ with }\quad   W_i := \left\{ x \in X \mid \lim_{z \to 0} z.x =x_i \right\},
\end{equation}
where $[m] := \{1,\dots,m\}$.  We call this a \textbf{BB-decomposition} since decompositions of this type were first studied by Bialynicki-Birula in \cite{Birula1973}. 

Now, we examine whether there exists a constructible grading such that the induced BB-decomposition is an affine paving. From \cite{LP23a}, we know that it has to satisfy the following necessary assumption.
\begin{dfn}\label{dfn:attractive-grading}
A grading $\mr{wt}:B \to \bZ$ on $Q(M,B)_0$ is \textbf{attractive} if:
\begin{itemize}
\item[(AG1)] for any $i\in Q_0$ it holds that $\mr{wt}(b^{(i)}_k)>\mr{wt}(b^{(i)}_\ell)$ whenever $k>\ell$,
\item[(AG2)] for any $a\in Q_1$ there exists a weight $\mr{wt}(a)\in\bZ$ such that \[\mr{wt}\big(b^{(t_a)}_\ell\big)=\mr{wt}\big(b^{(s_a)}_k\big)+\mr{wt}(a)\] whenever $b^{(s_a)}_k\to b^{(t_a)}_\ell\in Q(M,B)_1$.
\end{itemize}
\end{dfn} 
Clearly, this holds for the grading constructed in the proof of Lemma~\ref{lem:finitely-many-fixed-points} because we can reorder the basis so that it satisfies (AG1) and it already satisfies (AG2) by construction. Furthermore, the grading has to be compatible with the following condition on the coefficient quiver. This is necessary to prove that the describing equations of the BB-cells in the coordinates of $B$ parameterize an affine space.
\begin{dfn} A forest $M \in  \mr{rep}_\bC (Q)$ is \textbf{alignable} if there exists a basis $B$, such that for $Q(M,B)$ the following holds over each $(a : i \to j)\in Q_1$:
\begin{itemize}
\item[(SA1)] endpoints of segments have larger indices than points with outgoing arrows:
if $M_a b^{(i)}_\ell = 0$ and $M_a b^{(i)}_k \neq 0$, then $ \ell > k$. 
\item[(SA2)] outgoing arrows are order preserving:\\
if $M_a b^{(i)}_\ell =  b^{(j)}_{\ell'}$ and $M_a b^{(i)}_k = b^{(j)}_{k'}$ with $\ell>k$, then $ \ell' > k'$. 
\end{itemize}
\end{dfn}
\begin{prop}\label{prop:attractive-alignment}
    Let  $M \in \mathrm{rep}_\bC (Q)$ be straight. Then there exists a basis $B$ such that $Q(M,B)$ is aligned and there is a compatible attractive grading $\mr{wt}: B \to \bZ$.
\end{prop}
\begin{proof}
    We construct an alignment and attractive grading inductively. First, we place the end points of the segments over the vertices of $Q$, such that their index increases with increasing length of the corresponding segment. By (SA2) this order has to be preserved for all predecessors. Let $D$ be the maximal number of segments ending at a vertex. The endpoints are graded increasingly with increasing index by numbers in $[D] = \{1,\dots,D\}$ starting from one. To avoid conflicts in the grading we assume $\mr{wt}(a)=D$ for all $a \in Q_1$. Now we inductively add the other vertices going backwards along the equioriented strings. If at some step of this procedure two strings of the same length meet at the same vertex, their order is fixed based on a total order of the vertices $i \in Q_0$. Following the procedure described in the proof of Lemma~\ref{lem:finitely-many-fixed-points}, it is possible to shift the terminal weights of the segments to avoid conflicts in the grading. If the support of $M$ is quiver of type $\texttt{A}$ or a tree it is checked analogous to \cite[Proposition~5.1]{LP23a}, 
    that this is sufficient to produce an attractive grading. 

    If the support of $M$ is not a tree, $Q$ contains at least two parallel paths that have only the start and endpoint in common. Since we set $\mr{wt}(a)=D$ for all $a \in Q_1$, this essentially reduces to the type $\texttt{A}$ case if both paths are of the same length. If the parallel paths are of different length we have to reduce the edge weights along the longer path such that both paths have the same weight. Then this case also reduces to the type $\texttt{A}$ case studied in \cite[Proposition~5.1]{LP23a}.  

    The coefficient quiver constructed with the above procedure is aligned. If the endpoints of the segments satisfy (AG1), the inductive procedure based on (AG2) preserves that property, so that the final grading is attractive. 
\end{proof}
\begin{rem}
This grading describes a generic cocharacter of $T_M$.
\end{rem}
\begin{rem}
Since, the constructed grading depends on a fixed total order of $Q_0$ it is not unique. In fact, it is useful to work with different gradings to study GKM-variety structures (see Section~\ref{sec:KT-basis}).
\end{rem}
\begin{thm}\label{thm:Affine-Paving}
    Let $M \in \mathrm{rep}_\bC (Q)$ be straight and $\mb{e} \in \bZ^{Q_0}$. The action of $\bC^*$ as in Proposition~\ref{prop:attractive-alignment} induces an affine paving of $\mr{Gr}_\mb{e}(M)$.
\end{thm}
\begin{proof}
By \cite[Lemma~4.12]{Carrell02}, the BB-decomposition is an $\alpha$-partition. It remains to show that the attracting set of every fixed point is an affine cell. Let $V \in \mr{Gr}_\mb{e}(M)$ be $\bC^*$-fixed. By \cite[Theorem~1]{Cerulli2011} it is a coordinate subrepresentation, i.e. there exist subsets $K_i \in \binom{[m_i]}{e_i}$ such that $V^{(i)}$ is spanned by $\{ b^{(i)}_k : k \in K_i \}$. Here, $m_i :=\dim M^{(i)}$ and $\binom{[n]}{k}$ denotes the set of all $k$-element subsets of $\{1,\dots,n\}$. Let $a \in Q_1$ and $k\in[m_{s_a}]$ be such that $M_a b_k^{({s_a})}\neq 0$, then there exists an $\ell \in [m_{t_a}]$ such that $M_ab_k^{(s_a)}=b^{({t_a})}_\ell$. We use the notation $k':=\ell$ to stress the dependence on $k$.

Now, let $W$ be a point in the attracting set of $V$. We denote its basis vectors by $\{ w^{(i)}_k : k \in K_i \}$ for all $i \in Q_0$. By construction of the grading and definition of the  attracting sets, they are of the form 
\begin{equation}\label{eqn:CoordinateDescriptionPointCell}
 w_{k}^{(i)} = b_{k}^{(i)} + \sum_{j \in [m_i] \setminus  [k] \, : \ j \notin K_i} \mu_{j,k}^{(i)} b_j^{(i)}\qquad\hbox{with } \mu_{j,k}^{(i)} \in \bC.
\end{equation}
Analogous to the proof of \cite[Theorem~4.13]{Pue2022}, we compute that if $M_a b_{k}^{({s_a})} \neq 0$, the coefficients are subject to the conditions 
\begin{align}
\label{eqn1:CoordinateDescriptionConditions}\mu_{j,k}^{(s_a)} &= \mu_{j',k'}^{({t_a})} + \sum_{\substack{\ell \in [j-1] \setminus [k]\, :  \\ M_ab_\ell^{({s_a})} \neq 0, \\ \ell' \in K_{t_a} \setminus K_{s_a}' }} \mu_{\ell,k}^{({s_a})} \, \mu_{j',\ell'}^{({t_a})}  \quad \big( \mr{if} \ j \in [m_{s_a}]\setminus [k] \, : \ M_a b_{j}^{({s_a})}\neq 0, \ j' \notin K_{t_a}\big), \\
\label{eqn2:CoordinateDescriptionConditions} 0 &= \mu_{h,k'}^{({t_a})} + \sum_{\substack{\ell \in [m_{s_a}] \setminus  [k] \, : \\ M_ab_\ell^{({s_a})} \neq 0, \ \ell' < h, \\ \ell' \in K_{t_a} \setminus K_{s_a}', }} \mu_{\ell,k}^{({s_a})} \, \mu_{h,\ell'}^{({t_a})} \quad \binom{\mr{if} \ h \in [m_{t_a}] \setminus  [k'] \, : \ h \notin K_{t_a},}{\nexists \ell \in [m_{s_a}]\setminus[k] \ \mr{s.t.:}  \ M_a b_{\ell}^{({s_a})}= b_h^{({t_a})}}.
\end{align}
A parameter $\mu^{(i)}_{j,k}$ with $k \in K_i$ and $j \notin K_i$ is called initial if $b^{(i)}_k$ has no predecessor in $Q(M,B)$ or its predecessor is not a spanning vector of $V$, and for all successors of $b^{(i)}_j$ that do not span $V$ there exists a successor of $b^{(i)}_k$ over the same vertex. 

Since $Q(M,B)$ consists of equioriented strings, we can inductively resolve the above equations starting from the initial parameters. The initial parameters are independent since they have no predecessors. All remaining parameters have predecessors and hence depend on the initial parameters. Thus, the attracting set is an affine space and its dimension equals the number of initial parameters obtained from the $K_i$'s.
\end{proof}
\subsection{Torus orbits}
\begin{thm}\label{thm:one-dim-T-orbits}
    Let $M \in \mathrm{rep}_\bC (Q)$ be straight and $\mb{e} \in \bZ^{Q_0}$. The torus $T_M$ as in Definition~\ref{dfn:torus-action} acts on $\mr{Gr}_\mb{e}(M)$ with finitely many one-dimensional orbits.
\end{thm}
\begin{proof} 
We utilize the coordinate description of the attracting sets from the proof of Theorem~\ref{thm:Affine-Paving}. By construction, the attracting sets of the $\bC$-fixed points are $T_M$-stable. Moreover, the weight pairs of the vertices corresponding to the initial parameters are all distinct. Hence, the dimension of a $T_M$-orbit of a point in the attracting set equals its number of non-zero initial parameters. In particular, the number of one-dimensional $T_M$-orbits in the attracting set is equal to the number of initial parameters of the fixed point. Since the number of fixed points is finite, the total number of one-dimensional $T_M$-orbits is finite.
\end{proof}
\subsection{Proof of Theorem~\ref{thm:GKM-structure}}
   The action from Definition~\ref{dfn:torus-action} extends to the quiver Grassmannian by Proposition~\ref{prop:strong-admissible-gradings}. It is skeletal by Lemma~\ref{lem:finitely-many-fixed-points} and Theorem~\ref{thm:one-dim-T-orbits}. The affine paving from Theorem~\ref{thm:Affine-Paving} guarantees the vanishing of the odd-degree rational cohomology. Hence $(\mr{Gr}_\mb{e}(M),T_M)$ is a GKM-variety by Definition~\ref{dfn:GKM-variety}. 
\subsection{The non-straight setting}
Similar to Proposition~\ref{prop:flexible-reduction},  it is possible to collapse arrows $a:i\to j$, if $M_a=\mr{id}$ and $e_i=e_j$. We remove the vertex $j$, replace the arrows $b$ with $s_b=j$ by $b\circ a$ and replace $M_b$ by $M_b\circ M_a$. If this collapse is not possible, the pair $(\mb{e},M)$ is called \textbf{identity-free}.
\begin{lem}\label{lem:no-GKM-structure}
    Let $Q$ be a finite quiver, let $M \in \mr{rep}_\bC (Q)$ be a string representation and let $\mb{e} \in \bZ^{Q_0}$ be $M$-flexible. There exists no GKM-variety structure on $\mr{Gr}_\mb{e}(M)$ if $M$ is not straight. 
\end{lem}
\begin{proof}
    Since $M \in \mr{rep}_\bC (Q)$ is a string and $\mb{e} \in \bZ^{Q_0}$ is $M$-flexible, every $\mb{e}$-admissible grading is constructible by Lemma~\ref{lem:all-nice-gradings-are-constructible}. Without loss of generality, we can assume that the pair $(\mb{e},M)$ is identity-free. By the collapsing described above, this does not change the variety but reduces the number of cases we have to distinguish. 
    
    Now, if $M$ is not straight, the subquiver $ \bullet \rightarrow \bullet \leftarrow \bullet$ or $\bullet \leftarrow \bullet \rightarrow \bullet$ is contained in $\mr{supp}M \subseteq Q$. It supports at least three indecomposable summands of $M$ because $\mb{e}$ is $M$-flexible: The entries of the dimension vector $\mb{e}$ at these vertices are at least $(1,2,1)$ or $(2,1,2)$ and at most $\dim M^{(i)}-1$ for any $i \in  Q_0$. Observe that for the two-sink case the middle entry has to be strictly larger than the outer entries since $\mb{e} \in \bZ^{Q_0}$ is $M$-flexible and $(\mb{e},M)$ is identity-free. Analogously, in the two-source case, the middle entry is strictly smaller. The outer entries are allowed to differ in both cases. In the subsequent proposition we show that these two minimal cases do not admit GKM-variety structures. Hence, no GKM-structure can exist if they are contained.
\end{proof}
\begin{prop}
    Let $Q= 1 \rightarrow 2 \leftarrow 3$, $M = \mathbb{C}^3 \rightarrow \mathbb{C}^3 \leftarrow \mathbb{C}^3$ where every map is the identity and $\mathbf{e} = (1,2,1)$.  Then $\mathrm{Gr}_\mathbf{e}(M)$ does not admit a GKM-variety structure.
    \end{prop}
    
\begin{proof}     
Let $B = \{ b_{i,k} : i \in [3], k \in [3] \}$ be the basis of $M$ where $b_{i,k}$ denotes the $k$-th standard basis vector for the $i$-th copy of $\mathbb{C}^3$. By \cite[Theorem~1]{Cerulli2011}, the fixed points of the $\mathbb{C}^*$-action with $\lambda.b_{i,k} := \lambda^k \cdot b_{i,k}$ are given by coordinate subrepresentations. In particular,  
 \[ p_1 = \big(\langle b_{1,2}\rangle,\langle b_{2,2},b_{2,3}\rangle,\langle b_{3,2}\rangle\big) \]
is a $\mathbb{C}^*$-fixed point, where $\langle \dots \rangle$ denotes the span as $\mathbb{C}$-vector space. The attracting set of $p$ is 
\begin{align*}
C :=& \Big\{ q \in \mathrm{Gr}_\mathbf{e}(M) : \lim_{\lambda \to 0} \lambda.q = p\Big\}\\
=& \Big\{ q =  \big(\langle b_{1,2}+ab_{1,3}\rangle,\langle b_{2,2},b_{2,3}\rangle,\langle b_{3,2}+cb_{3,3}\rangle\big) : a,c \in \mathbb{C} \Big\}.
\end{align*}
From Corollary~\ref{cor:faithful-action} we obtain that the largest torus acting faithfully is of rank $3+2$. Without loss of generality we can assume that $t= (t_1,\dots,t_5) \in T = (\bC^*)^5$ acts as follows: 
\[ t.b_{1,k} := t_kt_4 \cdot b_{1,k}, \quad t.b_{2,k} := t_k \cdot b_{2,k}, \quad t.b_{3,k} := t_kt_5 \cdot b_{3,k},\]
For any $t\in T$ and any point $q \in C$, we compute 
\begin{align*}
t.q 
&= \big(\langle t_2t_4b_{1,2}+t_3t_4ab_{1,3}\rangle,\langle t_2b_{2,2},t_3b_{2,3}\rangle,\langle t_2t_5b_{3,2}+t_3t_5cb_{3,3}\rangle\big)\\
&= \big(\langle b_{1,2}+\frac{t_3}{t_2}ab_{1,3}\rangle,\langle b_{2,2},b_{2,3}\rangle,\langle b_{3,2}+\frac{t_3}{t_2}cb_{3,3}\rangle\big).
\end{align*}
Hence, the $T$-orbits of all points in $C\setminus\{p\}$ are one-dimensional. This implies that $C$ contains infinitely many one-dimensional $T$-orbits because $C$ is two-dimensional. Thus $\mathrm{Gr}_\mathbf{e}(M)$ does not admit a GKM-variety structure.
\end{proof}
\begin{rem}
For $M = \mathbb{C}^3 \leftarrow \mathbb{C}^3 \rightarrow \mathbb{C}^3$ with identity maps and $\mathbf{e} = (1,2,1)$ we show with the same arguments that $\mathrm{Gr}_\mathbf{e}(M)$ has no GKM-variety structure.
\end{rem}

Nevertheless, there is strong evidence that GKM-variety structures also exist for a suitable class of forests (see Appendix~\ref{app:example}). However, the construction of attractive gradings that are compatible with the aligned coefficient quiver and the study of the attracting sets require new methods if the coefficient quiver of $M$ contains trees. Hence, this is an independent topic and will not be covered in this paper.

\section{Combinatorial description of the moment graph and applications}\label{sec:momnent-graph}
In this section we provide a combinatorial description of the moment graphs corresponding to the GKM-varieties from Theorem~\ref{thm:GKM-structure}. This induces a description of the tangent spaces of these varieties and allows to describe a nice basis for their torus equivariant cohomology.
\subsection{Mutations of coefficient quivers}
A subquiver $S \subseteq Q$ is successor closed if for all $i \in S_0$ and all $a \in Q_1$ with $s_a=i$ also $t_a \in S_0$ holds. In particular, the successor closed subquiver $S$ is full, i.e. $S_1$ contains all edges of $Q_1$ connecting vertices in $S_0$. Hence the successor closed subquivers of a coefficient quiver $Q(M,B)$ are parameterized by a subset of $B$. Hence, they are in one to one correspondence with coordinate subrepresentations.  
\begin{dfn}\label{dfn:fundamental-mutation}
Let $M \in \mr{rep}_\bC (Q)$ be straight and $B$ a basis of $M$. A quiver morphism $m: S\to S'$ between two $\mb{e}$-dimensional successor closed subquivers $S,S' \subseteq Q(M,B)$ is called \textbf{mutation} if $L=S\setminus S'$ and $L'=S'\setminus S$ are connected and project to the same path in $Q$. The morphism $m$ is uniquely determined by $S$ together with the triple $(i,k,\ell)$, where $b^{(i)}_k$ is the starting vertex of $L$ and $b^{(i)}_\ell$ is the starting vertex of $L'$. We write $m^{(i)}_{k,\ell}$ and call it \textbf{fundamental mutation} if $k < \ell$.
\end{dfn}
\begin{rem}
  If $k>\ell$, the map $m^{(i)}_{k,\ell}$ is the inverse of a fundamental mutation. 
\end{rem}
\subsection{Moment graph description}We define the characters 
\[ \epsilon_q : T_M \to \bC^*, \quad \big((\gamma_i)_{i \in [d]},(\nu_j)_{j \in [c]}\big) \mapsto \gamma_q, \quad \text{ for } q \in[d],\]
\[ \delta_a : T_M \to \bC^*, \quad \big((\gamma_i)_{i \in [d]},(\nu_j)_{j \in [c]}\big) \mapsto \nu_a, \quad \text{ for } a \in[c].\]
\begin{thm}\label{thm:moment-graph}Let $M \in \mr{rep}_\bC (Q)$ be straight and $\mb{e} \in \bZ^{Q_0}$, and let $T_M$ act on $X=\mr{Gr}_\mb{e}(M)$ as in Definition~\ref{dfn:torus-action}. Take a generic attractive cocharacter $\chi\in \mathfrak{X}_*(T_M)$ and compatible $B$. 
The corresponding moment graph $\mathcal{G}=\mathcal{G}(X,T_M, \chi)$ is given by the following data:
\begin{itemize}
\item[(G0)] the vertices are $\mb{e}$-dimensional successor closed subquivers $S \subseteq Q(M,B)$,
\item[(G1)] the fundamental mutations $m^{(i)}_{k,\ell}:S\to S'$ parameterize $\mathcal{G}_1$,
\item[(G2)]  the edge $m^{(i)}_{k,\ell}:S\to S'$ is labeled by the character \[\epsilon_r+\delta_{a'_1}+\dots+\delta_{a'_{p'}}-(\epsilon_q+ \delta_{a_1}+\dots+\delta_{a_p}),\] where $r$ is the index of the segment containing $b^{(i)}_\ell$, $q$ is the index of the segment containing $b^{(i)}_k$, $a'_{p'}\circ \dots \circ a'_1$ is the path in $Q$ underlying the preceding segment of $b^{(i)}_\ell$ in $Q(M,B)$ and $a_p\circ \dots \circ a_1$ is the path in $Q$ underlying the preceding segment of $b^{(i)}_k$ in $Q(M,B)$.
\end{itemize}
\end{thm}
\begin{proof}
By Lemma~\ref{lem:finitely-many-fixed-points}, the $T_M$-fixed points are coordinate subrepresentations. The discussion above the theorem implies that coordinate subrepresentations are parameterized by $\mb{e}$-dimensional successor closed subquivers $S \subseteq Q(M,B)$. 

In Theorem~\ref{thm:one-dim-T-orbits} we prove that the one-dimensional $T_M$-orbits are parameterized by initial parameters $\mu^{(i)}_{\ell,k}$. These are in bijection with fundamental mutations $m^{(i)}_{k,\ell}$ because both definitions are based on the same structural analysis of subsets of $B$ respective subquivers of $Q(M,B)$.

By Definition~\ref{dfn:torus-action}, $T_M$ acts on $b^{(i)}_k$ with $\gamma_q\nu_{a_1}\cdot\ldots\cdot\nu_{a_p}$ and on $b^{(i)}_\ell$ it acts with $\gamma_r\nu_{a'_1}\cdot\ldots\cdot\nu_{a'_{p'}}$. Hence, we obtain the desired character by rescaling the coefficient of  $b^{(i)}_k$ to one.
\end{proof}
\subsection{Basis for the torus equivariant cohomology}\label{sec:KT-basis}
Let $(X,T)$ be a GKM-variety with generic cocharacter $\chi \in \mathfrak{X}_*(T)$. If there exists a directed path (possibly of length zero) from $x$ to $y$ in the moment graph $\mathcal{G}:=\mathcal{G}(X,T,\chi)$ write $x \succeq_\chi y$. This is a partial order if $\mathcal{G}$ is acyclic.
\begin{dfn}[{\cite[Definition 2.12]{Tymoczko2008}}]\label{def:K-T-class}
Let $(X,T)$ be a GKM-variety with acyclic moment graph $\mathcal{G}=\mathcal{G}(X,T,\chi)$. A \textbf{Knutson-Tao class} for $x \in X^T$ is an equivariant class $p^x = (p^x_y)_{y \in X^T}\in H_T^\bullet(X)$ such that:
\begin{enumerate}
\item[(KT1)] \(p^x_x=  \prod_{E\in \mathcal{G}_1: \ S_E=x} \alpha_E, \)
\item[(KT2)] each $p^x_y$ is a homogeneous polynomial in $\bQ[T]$ with $\deg p^x_y = \deg p^x_x$, 
\item[(KT3)] $p^x_y = 0$ for each $y \in X^T$ such that $y \nsucceq_\chi x$.
\end{enumerate}
\end{dfn}
A basis of $H_T^\bullet(X)$ consisting of KT-classes is called KT-basis.
\begin{lem}\label{lem:KT-basis}
   Let $M \in \mr{rep}_\bC (Q)$ be straight and $\mb{e} \in \bZ^{Q_0}$, and let $T_M$ act on $X=\mr{Gr}_\mb{e}(M)$ as in Definition~\ref{dfn:torus-action}. Take a generic attractive cocharacter $\chi\in \mathfrak{X}_*(T_M)$ and compatible $B$. Then, there exists a KT-basis of $H^\bullet_{T_M}(X)$.
\end{lem}
\begin{proof}
By Theorem~\ref{thm:moment-graph}, the edges of the moment graph are given by fundamental mutations. Hence, their definition implies that $\mathcal{G}$ contains no oriented cycles. The BB-decomposition is an affine paving by Theorem~\ref{thm:Affine-Paving}. Thus, cells are rational in the sense of \cite[Definition~1.9]{LP23a}, and \cite[Theorem~3.9]{LP23b} implies the existence of a KT-basis.
\end{proof}
A moment graph $\mathcal{G}$ is called Palais-Smale if for every $x \to y \in \mathcal{G}_1$, $x$ has more outgoing edges than $y$. It is shown in \cite[Lemma 2.16]{Tymoczko2008} that the KT-basis is unique if $\mathcal{G}$ is Palais-Smale. 

All known cases of quiver Grassmannians with this property are covered in \cite[Theorem~3.22]{LP23b}. In particular, the support of $M$ is required to be an equioriented quiver of type $\text{A}$ or $\tilde{\texttt{A}}$. So far it was not possible to find a suitable replacement for the assumptions of \cite[Theorem~3.22]{LP23b} in the context of straight representations. Hence, it would be very interesting to find one example which is not already covered by that statement. An example with a Palais-Smale moment graph where $M$ is not a string is given in Appendix~\ref{app:example}.
\subsection{Tangent space}\label{sec:tangent-space}
Following \cite{CaRe2008}, we obtain a description of the tangent spaces at a quiver Grassmannian by morphisms between quiver representations. Let $P$ be a set of paths in $\bC Q$ of length at least two such that all paths beyond a fixed length are in the ideal $I$ in the path algebra $\bC Q$ generated by $P$. The finite-dimensional algebra $\bC Q /I$ is called zero relation algebra (or bounded path algebra).
\begin{lem}\label{lma:tangent-space-formula}
Let $(Q,P)$ parametrize a zero relation algebra and let $M \in \mr{rep}_\bC(Q,P)$. For a subrepresentation $U$ of $M$ with dimension vector $\mb{e}$, the tangent space of the quiver Grassmannian $\mr{Gr}_\mb{e}(M)$ at $U$ is isomorphic to $\mathrm{Hom}_Q(U,M/U)$. 
\end{lem} 
\begin{rem}
The proof of this statement is analogous to the construction in \cite[Section~4.1]{CaRe2008}, where it is proved for representations of finite quivers without oriented cycles. 
For their arguments, it is sufficient to assume that the representation $M$ is finite-dimensional, and if the quiver contains oriented cycles, concatenations of maps $M_a$ along paths in $Q$ vanish for all paths exceeding a certain length. The generalized version of \cite[Lemma~2]{CaRe2008}, which is the first step for the proof of Lemma~\ref{lma:tangent-space-formula}, is proved in \cite[Proposition~2.5]{Pue2019}, and \cite[Proposition~6]{CaRe2008} generalizes in the same way.
\end{rem}
If $M \in \mr{rep}_\bC(Q)$ is straight, we can use mutations to describe the tangent space:
\begin{thm}\label{thm:tangent-space-from-mutations}
Let $M \in \mathrm{rep}_\bC (Q)$ be straight, $\mb{e} \in \bZ^{Q_0}$ and let the torus $T_M$ act on $\mr{Gr}_\mb{e}(M)$ as in Definition~\ref{dfn:torus-action}. A basis of the tangent space at the fixed point $U \in \mr{Gr}_\mb{e}(M)^{T_M}$ is given by the mutations of its coefficient quiver.
\end{thm}
The proof of the above theorem is based on the following description of morphisms between trees from 
\cite[Section~2]{C-B1989}:
\begin{lem}\label{lem:morph-of-trees}
Let $(Q,P)$ parametrize a zero-relation algebra and let $X= F_*V_S$ and $Y=E_*V_T$ be tree representations of $Q$. The space $\mathrm{Hom}_Q(X,Y)$ has as basis the set of triples $(S',\mu,T')$ satisfying the conditions:
\begin{itemize}

\item[(H1)] $S'$ and $T'$ are non-empty connected subquivers of $S$ and $T$, respectively.

\item[(H2)] $\mu: S' \to T'$ is a quiver isomorphism with $\mu \circ E = F$.

\item[(H3)] If an arrow $a \in S_1$ ends at a vertex $i \in S'_0$ then $a \in S'_1$.

\item[(H4)] If an arrow $b \in T_1$ starts at a vertex $i \in T'_0$ then $b \in T'_1$.
\end{itemize}
The homomorphism $f$ corresponding to the triple $(S',\mu,T')$ consists of the maps
\[ f_i : \bigoplus_{j \in F^{-1}(i)} \bC \to \bigoplus_{k \in E^{-1}(i)} \bC \]
for each $i \in Q_0$, where $(f_i)_{j,k} = 1$ if $j \in S'_0$ and $k = \mu(j)$, and otherwise $(f_i)_{j,k} = 0$.
\end{lem}
\begin{proof}[Proof of Theorem~\ref{thm:tangent-space-from-mutations}]
By Lemma~\ref{lma:tangent-space-formula}, the tangent space at $U$ is isomorphic to $\mathrm{Hom}_Q(U,M/U)$. Since $M$ is straight, the coordinate subrepresentation $U \in \mr{Gr}_\mb{e}(M)^{T_M}$ is straight and also $M/U$ is straight. In particular, they consist of trees. Hence, we can apply Lemma~\ref{lem:morph-of-trees} to compute a basis of $\mathrm{Hom}_Q(U,M/U)$. The combinatorial moves described by the conditions (H1),...,(H4) are in one to one correspondence with the mutations of coefficient quivers described in Definition~\ref{dfn:coeff-quiver}.  
\end{proof}
\begin{rem}
The description of the tangent space in Theorem~\ref{thm:tangent-space-from-mutations} applies to all quiver Grassmannians of forests with GKM-variety structure if their moment graph has a combinatorial description analogous to Theorem~\ref{thm:moment-graph}.    
\end{rem}
\subsection{Hall strata}
\begin{dfn}\label{dfn:hall-strata}
Let $Q$ be a finite quiver, $M \in \mr{Rep}_\bC(Q)$, $\mb{e} \in \bZ^{Q_0}$. The \textbf{Hall stratum} of $U \in \mr{Gr}_\mb{e}(M)$ is the set 
\[ \mathcal{S}_{U,M/U} :=\big\{ V \in \mr{Gr}_\mb{e}(M) : V \cong U \text{ and } M/V \cong M/U \big\}.\]
\end{dfn}
\begin{rem}
Hall strata are a refinement of the stratification from \cite[Lemma~2.4]{CFR12} where the condition on the quotient is omitted. If $(Q,P)$ parametrizes a zero-relation algebra and $M \in \mr{Rep}_\bC(Q,P)$ is injective, both strata coincide. Clearly, the tangent space as described in Lemma~\ref{lma:tangent-space-formula} is constant along the Hall strata. 
\end{rem}
\appendix
\section{Example}\label{app:example}
In this appendix, we provide an example of a quiver Grassmannian with the structure of a GKM-variety which is not covered by Theorem~\ref{thm:GKM-structure}. By \cite{CFR12}, the flag variety $\mathcal{F}l_4$ is isomorphic to the $\texttt{A}_3$ quiver Grassmannian for 
\[ M = \bC^4 \stackrel{\mr{id}}{\to} \bC^4 \stackrel{\mr{id}}{\to} \bC^4 \quad \text{and } \mb{e} = (1,2,3).\]
Theorem~\ref{thm:GKM-structure} implies that $T_M=(\bC^*)^{4+3}$ acts on $\mr{Gr}_\mb{e}(M)$ with finitely many fixed points. By Remark~\ref{rem:T-action-in-type-A}, we know that it is sufficient to set $\nu_j=1$ for all edges. From now on, $T_M$ denotes the corresponding subtorus of the torus from Definition~\ref{dfn:torus-action}. Hence, there is a bijection between $\mb{e}$-dimensional coordinate subrepresentations of $M$ and permutations in $\Sigma_4$: The $i$-th letter of the permutation (in one-line notation) equals the height of the segment starting over the $i$-th vertex of $\texttt{A}_3$. In particular,
\begin{center}
\begin{tikzpicture}[scale=.8]
\node at (-1.7,0.75) {$\leftrightarrow$};
\node at (-2.7,1) {$\mr{id}=$};
\node at (-2.7,0.5) {$(1234)$};

\draw[arrows={-angle 90}, shorten >=2.5]  (-1,0) -- (0,0); 
\draw[arrows={-angle 90}, shorten >=2.5]  (0,0) -- (1,0);
\draw[fill=black] (-1,0) circle(.08);
\draw[fill=black] (0,0) circle(.08);
\draw[fill=black] (1,0) circle(.08);

\draw[arrows={-angle 90}, shorten >=2.5]  (-1,0.5) -- (0,0.5); 
\draw[arrows={-angle 90}, shorten >=2.5]  (0,0.5) -- (1,0.5);
\draw[fill=white] (-1,0.5) circle(.08);
\draw[fill=black] (0,0.5) circle(.08);
\draw[fill=black] (1,0.5) circle(.08);

\draw[arrows={-angle 90}, shorten >=2.5]  (-1,1) -- (0,1); 
\draw[arrows={-angle 90}, shorten >=2.5]  (0,1) -- (1,1);
\draw[fill=white] (-1,1) circle(.08);
\draw[fill=white] (0,1) circle(.08);
\draw[fill=black] (1,1) circle(.08);

\draw[arrows={-angle 90}, shorten >=2.5]  (-1,1.5) -- (0,1.5); 
\draw[arrows={-angle 90}, shorten >=2.5]  (0,1.5) -- (1,1.5);
\draw[fill=white] (-1,1.5) circle(.08);
\draw[fill=white] (0,1.5) circle(.08);
\draw[fill=white] (1,1.5) circle(.08);
\end{tikzpicture} $\quad$
\begin{tikzpicture}[scale=.8]
\node at (-2.4,0.75) {$(3124) \leftrightarrow$};

\draw[arrows={-angle 90}, shorten >=2.5]  (-1,0) -- (0,0); 
\draw[arrows={-angle 90}, shorten >=2.5]  (0,0) -- (1,0);
\draw[fill=white] (-1,0) circle(.08);
\draw[fill=black] (0,0) circle(.08);
\draw[fill=black] (1,0) circle(.08);

\draw[arrows={-angle 90}, shorten >=2.5]  (-1,0.5) -- (0,0.5); 
\draw[arrows={-angle 90}, shorten >=2.5]  (0,0.5) -- (1,0.5);
\draw[fill=white] (-1,0.5) circle(.08);
\draw[fill=white] (0,0.5) circle(.08);
\draw[fill=black] (1,0.5) circle(.08);

\draw[arrows={-angle 90}, shorten >=2.5]  (-1,1) -- (0,1); 
\draw[arrows={-angle 90}, shorten >=2.5]  (0,1) -- (1,1);
\draw[fill=black] (-1,1) circle(.08);
\draw[fill=black] (0,1) circle(.08);
\draw[fill=black] (1,1) circle(.08);

\draw[arrows={-angle 90}, shorten >=2.5]  (-1,1.5) -- (0,1.5); 
\draw[arrows={-angle 90}, shorten >=2.5]  (0,1.5) -- (1,1.5);
\draw[fill=white] (-1,1.5) circle(.08);
\draw[fill=white] (0,1.5) circle(.08);
\draw[fill=white] (1,1.5) circle(.08);
\end{tikzpicture} $\quad$
\begin{tikzpicture}[scale=.8]
\node at (-2.4,0.75) {$(4321) \leftrightarrow$};

\draw[arrows={-angle 90}, shorten >=2.5]  (-1,0) -- (0,0); 
\draw[arrows={-angle 90}, shorten >=2.5]  (0,0) -- (1,0);
\draw[fill=white] (-1,0) circle(.08);
\draw[fill=white] (0,0) circle(.08);
\draw[fill=white] (1,0) circle(.08);

\draw[arrows={-angle 90}, shorten >=2.5]  (-1,0.5) -- (0,0.5); 
\draw[arrows={-angle 90}, shorten >=2.5]  (0,0.5) -- (1,0.5);
\draw[fill=white] (-1,0.5) circle(.08);
\draw[fill=white] (0,0.5) circle(.08);
\draw[fill=black] (1,0.5) circle(.08);

\draw[arrows={-angle 90}, shorten >=2.5]  (-1,1) -- (0,1); 
\draw[arrows={-angle 90}, shorten >=2.5]  (0,1) -- (1,1);
\draw[fill=white] (-1,1) circle(.08);
\draw[fill=black] (0,1) circle(.08);
\draw[fill=black] (1,1) circle(.08);

\draw[arrows={-angle 90}, shorten >=2.5]  (-1,1.5) -- (0,1.5); 
\draw[arrows={-angle 90}, shorten >=2.5]  (0,1.5) -- (1,1.5);
\draw[fill=black] (-1,1.5) circle(.08);
\draw[fill=black] (0,1.5) circle(.08);
\draw[fill=black] (1,1.5) circle(.08);
\end{tikzpicture} 
\end{center} 
Here, the vertices describing the subrepresentations are filled black. For this bijection, the cell dimension in the quiver Grassmannian equals the dimension of the corresponding Schubert cell. The Schubert variety equals the union of cells over all successors in the moment graph. 
Thus, every Schubert variety of $\mathcal{F}l_4$ has an induced GKM-variety structure from the restriction of the $T_M$-action on $\mr{Gr}_\mb{e}(M)$ (see \cite[Lemma~5.16]{LP23b}).
\begin{rem}
The above results apply to arbitrary Schubert varieties in $\mathcal{F}l_n$. 
\end{rem}
Let $X_{(3124)}$ denote the Schubert variety for the permutation $(3124)$. It is identified with the quiver Grassmannian $\mathrm{Gr}_{\mb{e}'}(M')$ for 
\[ \begin{tikzpicture}[scale=.8]
\node at (-2.1,0.5) {$Q' =$};

\draw[arrows={-angle 90}, shorten >=5, shorten <=5]  (-1,0) -- (0,0); 
\draw[arrows={-angle 90}, shorten >=5, shorten <=5]  (0,0) -- (1,0);
\draw[arrows={-angle 90}, shorten >=5, shorten <=5]  (0,0) -- (0,1);
\node at (-1,0) {$1$};
\node at (0,0) {$2$};
\node at (1,0) {$3$};
\node at (0,1) {$4$};

\node at (6,0.5) {$,\quad  M' = V_{Q'} \oplus V_{Q'} \oplus V_Q, \quad \mb{e}' = \begin{pmatrix}
     & 1& \\
     1 & 2 & 3
\end{pmatrix} $};
\end{tikzpicture} 
\]
Here $Q = 1 \to 2 \to 3$ and the representation $V_\bullet$ has $\bC$ over every vertex of the quiver and identity maps along every arrow. In fact, let $B'$ be a basis of $M$ such that 
\[ \begin{tikzpicture}[scale=.65]
\node at (-2.9,-0.75) {$Q(M',B') =$};

\draw[arrows={-angle 90}, shorten >=3, shorten <=3]  (-1,0) -- (0,0); 
\draw[arrows={-angle 90}, shorten >=3, shorten <=3]  (0,0) -- (1,0);
\draw[arrows={-angle 90}, shorten >=3, shorten <=3]  (0,0) -- (0.5,0.6);
\draw[fill=black] (-1,0) circle(.08);
\draw[fill=black] (0,0) circle(.08);
\draw[fill=black] (1,0) circle(.08);
\draw[fill=black] (0.5,0.6) circle(.08);

\draw[arrows={-angle 90}, shorten >=3, shorten <=3]  (-1,0-1) -- (0,0-1); 
\draw[arrows={-angle 90}, shorten >=3, shorten <=3]  (0,0-1) -- (1,0-1);
\draw[arrows={-angle 90}, shorten >=3, shorten <=3]  (0,0-1) -- (0.5,0.6-1);
\draw[fill=black] (-1,0-1) circle(.08);
\draw[fill=black] (0,0-1) circle(.08);
\draw[fill=black] (1,0-1) circle(.08);
\draw[fill=black] (0.5,0.6-1) circle(.08);

\draw[arrows={-angle 90}, shorten >=3, shorten <=3]  (-1,0-2) -- (0,0-2); 
\draw[arrows={-angle 90}, shorten >=3, shorten <=3]  (0,0-2) -- (1,0-2);
\draw[fill=black] (-1,0-2) circle(.08);
\draw[fill=black] (0,0-2) circle(.08);
\draw[fill=black] (1,0-2) circle(.08);
\end{tikzpicture} 
\]
and let $\lambda \in \bC^*$ and $(\gamma_1,\gamma_2,\gamma_3) \in (\bC^*)^3 =:T$ act as follows
\[  \begin{tikzpicture}[scale=.75]
\draw[arrows={-angle 90}, shorten >=6, shorten <=6]  (-1,0) -- (0,0); 
\draw[arrows={-angle 90}, shorten >=6, shorten <=6]  (0,0) -- (1,0);
\draw[arrows={-angle 90}, shorten >=6, shorten <=6]  (0,0) -- (0.5,0.6);
\node at (-1,0) {$\lambda$};
\node at (0,0) {$\lambda$};
\node at (1,0) {$\lambda$};
\node at (0.5,0.6) {$\lambda$};

\draw[arrows={-angle 90}, shorten >=6, shorten <=6]  (-1,0-1) -- (0,0-1); 
\draw[arrows={-angle 90}, shorten >=6, shorten <=6]  (0,0-1) -- (1,0-1);
\draw[arrows={-angle 90}, shorten >=6, shorten <=6]  (0,0-1) -- (0.5,0.6-1);
\node at (-1,0-1) {$\lambda^2$};
\node at (0,0-1) {$\lambda^2$};
\node at (1,0-1) {$\lambda^2$};
\node at (0.5,0.6-1) {$\lambda^2$};

\draw[arrows={-angle 90}, shorten >=6, shorten <=6]  (-1,0-2) -- (0,0-2); 
\draw[arrows={-angle 90}, shorten >=6, shorten <=6]  (0,0-2) -- (1,0-2);
\node at (-1,0-2) {$\lambda^3$};
\node at (0,0-2) {$\lambda^3$};
\node at (1,0-2) {$\lambda^3$};
\end{tikzpicture}, \quad 
 \begin{tikzpicture}[scale=.75]
\draw[arrows={-angle 90}, shorten >=6, shorten <=6]  (-1,0) -- (0,0); 
\draw[arrows={-angle 90}, shorten >=6, shorten <=6]  (0,0) -- (1,0);
\draw[arrows={-angle 90}, shorten >=6, shorten <=6]  (0,0) -- (0.5,0.6);
\node at (-1,0) {$\gamma_1$};
\node at (0,0) {$\gamma_1$};
\node at (1,0) {$\gamma_1$};
\node at (0.5,0.6) {$\gamma_1$};

\draw[arrows={-angle 90}, shorten >=6, shorten <=6]  (-1,0-1) -- (0,0-1); 
\draw[arrows={-angle 90}, shorten >=6, shorten <=6]  (0,0-1) -- (1,0-1);
\draw[arrows={-angle 90}, shorten >=6, shorten <=6]  (0,0-1) -- (0.5,0.6-1);
\node at (-1,0-1) {$\gamma_2$};
\node at (0,0-1) {$\gamma_2$};
\node at (1,0-1) {$\gamma_2$};
\node at (0.5,0.6-1) {$\gamma_2$};

\draw[arrows={-angle 90}, shorten >=6, shorten <=6]  (-1,0-2) -- (0,0-2); 
\draw[arrows={-angle 90}, shorten >=6, shorten <=6]  (0,0-2) -- (1,0-2);
\node at (-1,0-2) {$\gamma_3$};
\node at (0,0-2) {$\gamma_3$};
\node at (1,0-2) {$\gamma_3$};
\end{tikzpicture} 
\]
Then the induced GKM-variety structure of $T_M$ acting on $\mathrm{Gr}_{\mb{e}}(M)$ and the structure from the $T$-action on $\mathrm{Gr}_{\mb{e}'}(M')$ coincide. Additionally, the moment graph of the $T$-action on $\mathrm{Gr}_{\mb{e}'}(M')$ is Palais-Smale. This implies that the KT-basis from Lemma~\ref{lem:KT-basis} is unique.
\begin{rem}\label{rem:Schubert-varieties-as-quiver-Grass}
Conjecturally, all Schubert varieties can be identified with quiver Grassmannians for an enhanced type $\texttt{A}$ quiver. If the permutation avoids the pattern $(312)$, no enhancements are required. For every $(312)$ pattern in the permutation, we have to enhance $Q$ (the equioriented type $\texttt{A}$ quiver) at the vertex $i$ corresponding to the one in the pattern. In the representation $M = V_Q \oplus \dots \oplus V_Q$ we have to enhance all segments above the subsegment of the permutation starting at the $i$-th vertex. The new entries of the dimension vector $\mb{e}'$ count the enhanced subsegments corresponding to the permutation. The $\bC^*$- and $T$-action are analogous to the above example and induce a GKM-variety structure.
\end{rem}
\begin{rem}
Following Remark~\ref{rem:Schubert-varieties-as-quiver-Grass}, the above example is one instance of a class of quiver Grassmannians with GKM-varity structure where the representations are no strings. It would be very interesting to classify these representations. One possible subcase is the restriction to representations supported on a tree. 
\end{rem}
\section*{Acknowledgements} 
The section on constructible gradings is inspired by discussions with Martina Lanini and Giovanni Cerulli Irelli. I would like to thank Alexander Samokhin for explaining the connection between full exceptional collections and moment graphs of GKM-varieties. The example in the Appendix was developed in discussions with Giulia Iezzi. I was funded by the Deutsche Forschungsgemeinschaft (DFG, German Research Foundation) — SFB-TRR 358/1 2023 — 491392403. 

\end{document}